\documentclass [12pt,a4paper]{article}
\usepackage {amsmath,amscd,amsfonts,latexsym,theorem}

\newcommand {\cc}{\mathcal {C}}
\newcommand {\OO}{\mathcal {O}}
\newcommand {\CC}{\mathbb {C}}
\newcommand {\HH}{\mathbb {H}}
\newcommand {\PP}{\mathbb {P}}
\newcommand {\ZZ}{\mathbb {Z}}
\newcommand {\Hom}{\operatorname {Hom}}
\newcommand {\Pic}{\operatorname {Pic}}

\newcommand {\EExt}{\operatorname {\mathcal{E}\mathit{xt}}}
\newcommand {\coker}{\operatorname {coker}}
\newcommand {\im}{\operatorname {im}}
\newcommand {\codim}{\operatorname {codim}}
\newcommand {\vdim}{\operatorname {virt~dim}}
\newcommand {\mnx}{\bar M_N (X,\beta)}
\newcommand {\mnp}{\bar M_N (\PP^n,dH')}
\newcommand {\mnxx}[1]{\bar M_N (X,#1)}
\newcommand {\mnb}{\bar M_N (\tilde X,\beta-eE')}
\newcommand {\mnpb}{\bar M_N (\tilde \PP^n,dH'-eE')}
\newcommand {\mnbb}[1]{\bar M_N (\tilde X,#1)}

\newcommand {\mtx}{M_\tau (X,\beta)}
\newcommand {\mtb}{M_\tau (\tilde X,\beta-eE')}

\newcommand {\gwinv}[3]{\Phi_{#1}(#2 \; | \; #3)}

\newcommand {\ra}{\begin {picture}(15,0)
  \put (2,1){\line (1,0){10}}
  \put (10.5,-0.05){$\rightarrow$}
  \end {picture}}
\newcommand {\da}{\begin {picture}(0,10)
  \put (-3.5,0.5){\line (-1,1){7.7}}
  \put (3.5,0.5){\line (1,1){7.5}}
  \put (-7.04,0.9){$\searrow$}
  \put (8.3,6.25){$\nearrow$}
  \end {picture}}

\newtheorem {prop}{Proposition}[section]
\newtheorem {lemma}[prop]{Lemma}
{ \theorembodyfont {\rmfamily}
  \newtheorem {definition}[prop]{Definition}
}
\newenvironment {proof}{\textsc {Proof.~}}{~$\Box$\par\vspace{4mm}}
\newenvironment {example}{\textsc {Example.~}}{\par\vspace{4mm}}
\newenvironment {remark}{\textsc {Remark.~}}{\par\vspace{4mm}}

\begin {document}

%%%%%%%%%%%%%%%%%%%%%%%%%%%%%%%%%%%%%%%%%%%%%%%%%%%%%%%%%%%%%%%%%%%%%%%%%%%%%%%

\centerline {\Huge Counting rational curves with multiple points}
\vspace {2mm}
\centerline {\Huge and Gromov-Witten invariants of blow-ups}
\vspace {8mm}
\centerline {\large A. Gathmann}
\vspace {5mm}
\centerline {University of Hannover, Germany}
\centerline {e-mail: gathmann\symbol{64}math.uni-hannover.de}
\vspace {7mm}
\centerline {September 13, 1996}
\vspace {8mm}

\begin {abstract}
  We study Gromov-Witten invariants on the blow-up of $ \PP^n $ at a point,
  which is probably the simplest example of a variety whose moduli spaces of
  stable maps do not have the expected dimension. It is shown that many of
  these invariants can be interpreted geometrically on $ \PP^n $ as certain
  numbers of rational curves having a multiple point of given order at the
  blown up point. Moreover, all these invariants can actually be calculated,
  giving enumerative invariants of $ \PP^n $ which have not been known before. 
\end {abstract}

\vspace {5mm}

%%%%%%%%%%%%%%%%%%%%%%%%%%%%%%%%%%%%%%%%%%%%%%%%%%%%%%%%%%%%%%%%%%%%%%%%%%%%%%%

\section {Introduction} \label {intro}

Over the last few years, Gromov-Witten invariants have become a very powerful
tool in enumerative geometry. Let us briefly recall their definition. If $X$
is an $n$-dimensional smooth complex projective variety and $ \beta \in H_2
(X) $ a homology class, then for every $ N \ge 3 $ one defines a moduli space
$ \mnx $ of genus zero $N$-pointed stable maps into $X$, which is a
compactification of the space of all maps $ f: C \to X $, where $C$ is a
smooth $N$-pointed rational curve (which is allowed to vary) \cite {BM}. If
$X$ is a convex variety, i.e.\ if $ h^1 (C, f^* T_X) = 0 $ for all $ f: C \to
X $, then this moduli space is a smooth Deligne-Mumford stack of the expected
dimension
  \[ \vdim \mnx = - K_X \cdot \beta + n - 3 + N. \]
Now if $ \alpha_1,\dots,\alpha_N \in A^*(X) $ are cohomology classes whose
codimensions sum up to $ \vdim \mnx $, one defines the associated
Gromov-Witten invariant by
  \[ \gwinv {X}{\alpha_1,\dots,\alpha_N}{\beta} =
       (p_1^* \alpha_1 \cdots p_N^* \alpha_N) \cdot [\mnx] \qquad (*) \]
where $ p_i: \mnx \to X $ are the obvious evaluation maps and $ [\mnx] $
denotes the fundamental class of $ \mnx $. Geometrically, this invariant can
be interpreted as the number of rational curves in $X$ of homology class
$ \beta $ which pass through generically chosen subvarieties $ V_i \subset X $
with $ [ V_i ] = \alpha_i $.

If $X$ is not a convex variety, however, the actual dimension of $ \mnx $ will
in general be greater than the expected one, such that the above definition of
the Gromov-Witten invariants is not applicable. In this case, K. Behrend \cite
{B} has shown recently that it is possible to define a virtual fundamental
class
  \[ [\mnx]^{virt} \in A_{\vdim \mnx} (\mnx) \]
such that, if one uses this class in $(*)$ instead of the usual fundamental
class, this defines Gromov-Witten invariants (satisfying the usual axioms
\cite {KM}) on an arbitrary smooth complex projective variety $X$. Of course,
in this case there is no longer an obvious geometric interpretation of the
invariants.

In this paper, we study this construction in the case where $ X = \tilde \PP^n
$ is the blow-up of projective $n$-space in a point $ P \in \PP^n $. Let $ H' $
be the class of a line in $ \PP^n $ and $ E' $ be the class of a line in the
exceptional divisor $E$. We consider the commutative diagrams
  \[ \begin {CD}
         \mnpb       @>{\phi}>>    \mnp   \\
    @V{\tilde p_i}VV            @VV{p_i}V \\
       \tilde \PP^n  @>{\pi}>>     \PP^n
  \end {CD} \]
and show that, although there are components in $ \mnpb $ whose dimension is
too large, these are actually mapped by $ \phi $ to a subspace in $ \mnp $
whose dimension is \emph {smaller} than the expected one, which means that
they are irrelevant if one can compute the intersection products on the moduli
space $ \mnp $. This is obviously the case if all the classes $ \alpha_i $
in the Gromov-Witten invariant are pullbacks of classes on $ \PP^n $. Hence,
in this case it will again be possible to give a geometric interpretation of
the invariants.

It is even possible to give an interpretation of the Gromov-Witten invariants
on $ \tilde \PP^n $ in terms of curves on $ \PP^n $: via strict transform,
curves of degree $d$ in $ \PP^n $ which pass through $P$ with total
multiplicity $e$ correspond to curves in $ \tilde \PP^n $ of homology class
$ dH'-eE' $. This will lead to our main result (proposition \ref {meaning}):

\begin {it}
  Let $ d>0 $, $ e \ge 0 $ and $ \alpha_1,\dots,\alpha_N \in A^{\ge 1}
  (\PP^n) $ such that
    \[ \sum_i (\codim \alpha_i - 1 ) = d \, (n+1) - e \, (n-1) + n - 3. \]  
  Let $ P \in \PP^n $ be a point. Choose generic subschemes $ V_i \subset
  \PP^n $ with $ [V_i] = \alpha_i $ (in a sense that will be made precise).
             
  Then the number of rational curves in $ \PP^n $ (purely 1-dimensional
  subvarieties birational to $ \PP^1 $) of degree $d$ which have non-empty
  intersection with all $ V_i $ and which pass through the point $P$ with
  total multiplicity $e$, where each such curve $C$ is counted with
  multiplicity
    \[ \sharp (C \cap V_1) \cdots \sharp (C \cap V_N), \]
  is equal to the Gromov-Witten invariant on $ \tilde \PP^n $
    \[ \gwinv {\tilde \PP^n}{\pi^*\alpha_1,\dots,\pi^*\alpha_N}{dH' - eE'}. \]
\end {it}

Moreover, it will be shown that all Gromov-Witten invariants of $ \PP^n $
can be computed using the First Reconstruction Theorem of Kontsevich and Manin
\cite {KM} and some initial data that we will calculate.

In the case $ n=e=2 $, our results reproduce the numbers of rational curves
of degree $d$ in $ \PP^2 $ having a node in $P$ and passing through $ 3d-3 $
further points in the plane, which have already been computed last year by
R. Pandharipande \cite {P} using different methods.

I have been informed that L. G\"ottsche and R. Pandharipande have been
working on Gromov-Witten invariants of multiple blow-ups of $ \PP^2 $ and
their geometric interpretation.

The paper is organized as follows: In section \ref {map-image} we will give a
correspondence between stable maps to $X$ and embedded curves in $X$. In
section \ref {obstruction} we calculate the cohomology $ h^1 (C, f^* T_{\tilde
\PP^n}) $ for all maps $ f: C \to \tilde \PP^n $ which will enable us in
section \ref {image-map} to prove the statement on the dimension of the image
of the map $ \phi: \mnpb \to \mnp $ that was mentioned above. The proof of
our main result on the geometric meaning of the invariants on $ \tilde \PP^n $
will be given in section \ref {geometric}. Finally, in section \ref
{calculation} we show how to calculate the invariants on $ \tilde \PP^n $ and
give some examples.

\textbf {Notations and Conventions.} We will always work over the ground field
$ \CC $ of complex numbers. If $X$ is an $n$-dimensional smooth projective
variety, we follow \cite {BM} and let
  \[ H_2^+(X) = \{ \beta \in \Hom_{\ZZ} (\Pic X,\ZZ) \; | \;
       \beta(L) \ge 0 \;\, \mbox {whenever $L$ is ample} \} \]
be the semigroup of "positive homology classes" in $X$. If $ N>3 $ and $ \beta
\in H_2^+ (X) $, we denote by $ \mnx $ the Deligne-Mumford stack of $N$-pointed
stable maps in $X$ of genus zero and homology class $ \beta $ as defined in
\cite {BM}. In this paper, the terminology \emph {stable map} will always be
used to denote a stable map of genus zero. A stable map $ (C,x_1,\dots,x_N,f)
$ will be called \emph {irreducible} if the curve $C$ is, and reducible
otherwise.

If $ \alpha \in A^c(X) $ is a cycle, we denote the \emph {codimension} of $
\alpha $ by $ c = \codim \alpha $. If $ \alpha_1,\dots,\alpha_N \in A^* (X) $
are cycles on $X$ whose codimensions sum up to the virtual dimension
  \[ \vdim \mnx = - K_X \cdot \beta + n - 3 + N \]
of $ \mnx $, then K. Behrend \cite {B} has defined an associated Gromov-Witten
invariant which is denoted by
  \[ \gwinv {X}{\alpha_1,\dots,\alpha_N}{\beta}
       = (p_1^*\alpha_1 \cdots p_N^*\alpha_N) \cdot [\mnx]^{virt} \]
where $ p_i : \mnx \to X $ are the evaluation maps and $ [\mnx]^{virt} $
is the virtual fundamental class of the moduli space $ \mnx $.

%%%%%%%%%%%%%%%%%%%%%%%%%%%%%%%%%%%%%%%%%%%%%%%%%%%%%%%%%%%%%%%%%%%%%%%%%%%%%%%

\section {Stable maps and their images} \label {map-image}

Gromov-Witten invariants are concerned with stable maps into $X$. Since our
final aim is to make statements about the numbers of rational curves embedded
in a variety $X$, we start by collecting some relations between these two
points of view. So let us begin by defining the moduli spaces between which we
will find a correspondence later.

\begin {definition} \label {relevant}
  Let $X$ be a smooth projective variety over $\CC$, $ 0 \neq \beta \in H_2^+
  (X) $ a fixed homology class, and $ N \ge 3 $.
  \begin {enumerate}
    \item A \emph {relevant map} in $X$ is a (genus zero) stable map $
      (C,x_1,\dots,x_N,f) $ such that $ [f(C)] = f_*[C] = \beta \in H_2^+(X) $,
      i.e.\ there is no irreducible component of $C$ on which $f$ is a finite
      covering. (The map $f$ may, however, contract some irreducible components
      of $C$ to a point.) The set of all such relevant maps in $X$ modulo
      isomorphism can be regarded as a substack of $ \mnx $. It will be denoted
      by $ RM_N (X,\beta) $.
    \item The set of all irreducible relevant maps (i.e.\ relevant maps as
      above whose underlying curve $C$ is irreducible) such that $
      f^{-1}(f(x_i)) = \{ x_i \} $ for all $i$ will be denoted by $ RM'_N
      (X,\beta) $.
    \item A \emph {relevant curve} in $X$ is defined to be a tuple $ (D,y_1,
      \dots,y_N) $, where $D$ is a (purely) 1-dimensional, closed, connected
      subvariety of $X$ with $ [D] = \beta $ such that every irreducible
      component of $D$ is rational (i.e.\ birational to $ \PP^1 $, not
      necessarily smooth), and where the $ y_i $ are points on $D$, not
      necessarily distinct. Let $ RC_N (X,\beta) $ be the set of all such
      relevant curves.
  \end {enumerate}
\end {definition}

\begin {remark}
  We call such maps "relevant maps" because we will show later (see
  proposition \ref {finite}) that, under favourable circumstances, we can
  arrange that all curves counted by the Gromov-Witten invariants are of this
  type.
\end {remark}

The first property we want to show is that every relevant curve is the image
of some relevant map. This and the following results in this section will be
set-theoretic, since this is all we need.

\begin {prop} \label {rel-1}
  Let $X$ be a smooth projective variety, $ 0 \neq \beta \in H_2^+(X) $, and $
  N \ge 3 $. Then there is a natural surjective map
    \[ \mu: RM_N (X,\beta) \to RC_N (X,\beta) \]
  which is given by
    \[ \mu (C,x_1,\dots,x_N) = (f(C),f(x_1),\dots,f(x_N)). \]
\end {prop}

\begin {proof}
  Let $ (C,x_1,\dots,x_N,f) $ be a relevant map in $X$. We have to show that
  $ (f(C),f(x_1),\dots,f(x_N)) $ is a relevant curve.
  \begin {itemize}
    \item Obviously, $ f(C) \subset X $ is a 1-dimensional, closed, connected
      subvariety of $X$ with $ f(x_i) \in f(C) $.
    \item By definition of a relevant map, $ [f(C)] = f_* [C] = \beta $.
    \item Every irreducible component of $ f(C) $ can be written as $ f(C_0) $
      with $ C_0 $ a rational irreducible component of $C$. Since $ C_0 $
      is rational, so is $ f(C_0) $.
  \end {itemize}
  This shows that $ (f(C),f(x_1),\dots,f(x_N)) $ is a relevant curve.

  To show the surjectivity of the map $\mu$, let $ (D,y_1,\dots,y_N) $ be a
  relevant curve in $X$ and $ D = D_1 \cup \cdots \cup D_m $ its decomposition
  into irreducible components. 

  By assumption, each $ D_i $ is a rational curve, hence we can find rational
  maps $ f_i: C_i \to D_i $ with $ C_i \cong \PP^1 $ for all $i$. Of course,
  the $ f_i $ have to be surjective morphisms. Moreover, we have $ f_* [C_i] =
  D_i $, since the $ f_i $ are birational.

  Now the procedure to glue the $ f_i: C_i \to X $ to a map $ f: C \to X $
  and to choose points $ x_1,\dots,x_N \in C $ such that $ (C,x_1,\dots,x_N,f)
  $ becomes a relevant map with image $ (D,y_1,\dots,y_N) $ under $\mu$ is
  rather obvious, but nevertheless we will give it in detail.

  Since $D$ is connected, we can assume that its components $ D_i $ are
  numbered in such a way that for each $ i = 2,\dots,m $ there is an $ a(i)
  \in \{ 1,\dots,i-1 \} $ such that there exists a point $ y'_i \in D_i
  \cap D_{a(i)} $.

  Now, for every $ y \in D $ which is either one of the $ y_k $ or one of the
  $ y'_i $, do the following:

  If $ y=y'_i $ for exactly one $ i = 2,\dots,m $, but $ y \neq y_k $ for all
  $ k = 1,\dots,N $, then glue $ C_i $ and $ C_{a(i)} $ together at a point
  which is mapped to $y$ both by $ f_i $ and by $ f_{a(i)} $.

  If $ y=y_k $ for exactly one $ k = 1,\dots,N $, but $ y \neq y'_i $ for all
  $ i = 2,\dots,m $, then choose some point $ x_k $ in some $ C_i $ with
  $ f_i(x_k) = y $.

  In all other cases, let $ C'_y $ be a smooth rational curve and let $ f'_y:
  C'_y \to \{y\} $ be the constant map. If $ y \neq y'_i $ for all $ i = 2,
  \dots,m $, glue $ C'_y $ to some point in some $ C_i $ which is mapped to
  $y$. Otherwise, for any component $ C_i $ such that $ y=y'_i $ or $ i = a(j)
  $ for some $j$ with $ y=y'_j $, choose a point on this component which is
  mapped to $y$ and glue it at this point to some point on $ C'_y $. In both
  cases, for every $k$ with $ y=y_k $ choose a point $ x_k \in C'_y $. (All
  points chosen on $ C'_y $, the $ x_k $ as well as the points glued to the
  $ C_i $, have to be distinct.)

  Now let $C$ be the union of all $ C_i $ and $ C'_y $, glued together as
  described above. Let $ f: C \to X $ be the map induced by $ f_i $ and $ f'_y
  $. Then, as can be seen from the construction, $ (C,x_1,\dots,x_N,f) $ is a
  relevant map with image $ (D,y_1,\dots,y_N) $ under $\mu$. Hence, the map
  $\mu$ is surjective.
\end {proof}

The next result will eventually allow us to establish a one-to-one
correspondence between certain stable maps and their images.

\begin {lemma} \label {rel-2}
  Let 
    \[ \mu' = \mu|_{RM'_N (X,\beta)} : RM'_N (X,\beta) \to RC_N (X,\beta) \]
  be the restriction of the map considered in proposition \ref {rel-1}. Then
  \begin {enumerate}
    \item the image of $ \mu' $ is contained in the subset of $ RC_N (X,\beta)
      $ parametrizing irreducible curves.
    \item $ \mu' $ is injective.
    \item The stable maps in $ RM'_N (X,\beta) $ have no non-trivial
      automorphisms.
  \end {enumerate}
\end {lemma}

\begin {proof}
  If $ (C,x_1,\dots,x_N,f) \in RM'_N (X,\beta) $, then $C$ is irreducible,
  so $ f(C) $ is irreducible, too. This shows (i).

  Let $ (C,x_1,\dots,x_N,f) $ and $ (C',x'_1,\dots,x'_N,f') $ be two relevant
  maps in $ RM'_N (X,\beta) $ having the same image
    \[ (f(C),f(x_1),\dots,f(x_N)) = (f'(C'),f'(x'_1),\dots,f'(x'_N)). \]
  Since $ C \cong C' \cong \PP^1 $ and $f$, $f'$ are birational, we can
  consider the composition $ \alpha: f^{-1} \circ f' $ which is a also
  a birational map and therefore an isomorphism between $ C' $ and $C$.
  Because of the condition $ f^{-1}(f(x_i)) = \{ x_i \} $, $ \alpha $ maps
  each $ x'_i $ to $ x_i $ and hence induces an isomorphism between the two
  relevant maps, so they represent the same element in $ RM'_N (X,\beta) $,
  which proves (ii). Finally, the isomorphism constructed is obviously
  unique, which shows (iii).
\end {proof}

%%%%%%%%%%%%%%%%%%%%%%%%%%%%%%%%%%%%%%%%%%%%%%%%%%%%%%%%%%%%%%%%%%%%%%%%%%%%%%%

\section {Calculation of the obstruction $ h^1 (C,f^* T_{\tilde X}) $}
  \label {obstruction}

We now start to study the relation between stable maps on a variety $X$ and
on its blow-up $ \tilde X $ in a point $ P \in X $. Our main problem will be
that $ \tilde X $ is never a convex variety in the sense of \cite {K}, since
there are always stable maps $ (C,x_1,\dots,x_N,f) $ whose image is contained
in the excepional divisor where $ h^1 (C, f^* T_{\tilde X}) $ does not vanish
(see e.g.\ lemma \ref {coh-1}). Hence, we expect the moduli spaces of stable
maps into $ \tilde X $ to have the "wrong" dimension. Indeed, it is easy to
give an example, even for curves that do not lie entirely in the exceptional
divisor:

\begin {example}
  Let $ X = \PP^n $ and denote by $H'$ the class of a line in $X$. We consider
  curves of degree $d$, whose moduli space is well known to have the expected
  dimension, namely
    \[ \dim \mnxx {dH'} = \vdim \mnxx {dH'} = d \, (n+1) + n - 3 + N. \]
  Now consider curves on the blow-up $ \tilde X $ having homology class $ dH'
  $. Of course, the virtual dimension remains the same as above, but now we
  have new possibilities of realizing curves with this homology class using
  reducible curves where parts of it are lying in the exceptional divisor $E$:
  for example, take the strict transform $ C_1 $ of any curve of degree $d$ on
  $X$ passing $e$ times through the blown up point. This curve has homology
  class $ dH'-eE' $, where $E'$ denotes the class of a line in $ E \cong
  \PP^{n-1} $. Since passing through $P$ gives $ n-1 $ conditions, the
  dimension of the space of such curves is (at least)
    \[ D_1 = d \, (n+1) - e \, (n-1) + n - 3. \]
  Now take any curve $ C_2 $ in the exceptional divisor $ E \cong \PP^{n-1} $
  of degree $e$, the space of such curves has dimension $ en + n - 4 $. To be
  able to glue it to $ C_1 $, we have to require it to pass through one of the
  points where $ C_1 $ meets the exceptional divisor, which gives $ n-2 $
  conditions on $ C_2 $. Hence, for the choice of $ C_2 $, we have
    \[ D_2 = en - 2 \]
  degrees of freedom. Now we can consider reducible curves made up of two
  components $ C_1 $ and $ C_2 $ as above, these curves have homology class
  $ dH' $ and hence give, together with $N$ marked points, a subspace of
  $ \mnbb {dH'} $. But the dimension of this space is (at least)
    \[ D_1 + D_2 + N = d \, (n+1) + n - 3 + N + e - 2 \]
  which is bigger than the expected dimension of $ \mnbb {dH'} $ if $ e>2 $.
\end {example}

It will be the aim of this section to calculate $ h^1 (C,f^* T_{\tilde X}) $
if $ f: C \to \tilde X $ is a map from a prestable curve of genus zero to $
\tilde X $. This will tell us more precisely which parts of the moduli spaces
have the "correct" dimension.

First we introduce some notation which will be used throughout the rest of the
paper when dealing with blow-ups. Let $X$ be a smooth $n$-dimensional projective
variety, $ n \ge 2 $. We will soon specialize to the case where $ X = \PP^n $,
but for the moment we keep it arbitrary.

Let $ \tilde X $ be the blow-up of $X$ in a fixed point $ P \in X $, such that
we have a cartesian diagram
\[ \begin {CD}
   E   @>{i}>> \tilde X  \\
  @VVV         @VV{\pi}V \\
   P    @>>>       X
\end {CD} \]
By a Mayer-Vietoris type argument (see e.g.\ \cite {GH}, ch.\ 4.1) one sees
that
  \[ H_2 (\tilde X) = H_2 (X) \oplus \ZZ \cdot E' \]
where $ E' $ is the class of a line in the exceptional divisor $ E \cong
\PP^{n-1} $.

In the case $ X = \PP^n $ we will denote both the hyperplane class of $X$ as
well as its pullback to $ \tilde X $ by $H$. The class of a line will be
denoted $ H' $.

We now start by proving a lemma which will allow us in favourable cases to
reduce the calculation of $ h^1 (C,f^* T_{\tilde X}) $ to computations that
do not involve $ \tilde X $ but only $X$.

\begin {lemma} \label {triangle}
  Let $C$ be a smooth curve and $ f: C \to \tilde X $ a morphism such that
  $ f(C) \not\subset E $. Let $ D := f^* (E) $, which is an effective divisor
  on $C$. Then there is a commutative diagram of sheaves on $C$
  \[ \begin {array}{c}
    f^* \pi^* T_X (-D) \ra f^* \pi^* T_X \\
    ~~~~~~ \da \\
    ~~~~~~ f^* T_{\tilde X}
  \end {array} \]
  where all the three morphisms are injective, and all of them are
  isomorphisms away from the support of $D$.
\end {lemma}

\begin {proof}
  Since $ E = P \times_X \tilde X $, we have $ i^* \Omega_{\tilde X/X} =
  \Omega_{E/P} = \Omega_E $. As $ \Omega_{\tilde X/X} $ has support on $E$,
  this can be rewritten as $ i_* \Omega_E = \Omega_{\tilde X/X} $. Hence,
  there is an exact sequence of sheaves on $ \tilde X $
    \[ 0 \to \pi^* \Omega_X \to \Omega_{\tilde X} \to i_* \Omega_E \to 0. \]
  Dualizing, we get
    \[ 0 \to T_{\tilde X}
         \to \pi^* T_X
         \to \EExt^1 (i_* \Omega_E, \OO_{\tilde X})
         \to 0. \]
  By duality (\cite {HR}, thm.\ III 6.7), we have
    \[ \EExt^1 (i_* \Omega_E, \OO_{\tilde X})
         = i_* \EExt^1 (\Omega_E, N_{E/\tilde X})
         = i_* T_E (-1), \]
  therefore we get a morphism $ \pi^* T_X \to i_* T_E (-1) $ which we can
  restrict to $E$ to get a morphism $ \pi^* T_X|_E \to i_* T_E (-1) $ fitting
  into a commutative diagram
  \[ \begin {CD}
    0 @>>> \pi^* T_X (-E) @>>> \pi^* T_X @>>> \pi^* T_X|_E @>>> 0  \\
    @.           @.                 @|            @VVV          @. \\
    0 @>>>  T_{\tilde X}  @>>> \pi^* T_X @>>> i_* T_E (-1) @>>> 0
  \end {CD} \]
  From this we can deduce the existence of a map $ \pi^* T_X (-E) \to
  T_{\tilde X} $ giving a commutative diagram
  \[ \begin {array}{c}
    \pi^* T_X (-E) \ra \pi^* T_X \\
    ~~~~~~ \da \\
    ~~~~~~ T_{\tilde X}
  \end {array} \]
  with all three morphisms injective. Finally, apply the functor $ f^* $ to
  get the desired result
  \[ \begin {array}{c}
    f^* \pi^* T_X (-D) \ra f^* \pi^* T_X \\
    ~~~~~~ \da \\
    ~~~~~~ f^* T_{\tilde X}
  \end {array} \]
  Since $ f(C) \not\subset E $ by assumption, the morphisms are still
  injective, and certainly they are also isomorphisms away from the support of
  $D$.
\end {proof}

We are now ready to compute the relevant cohomology groups in the case where
$ C = \PP^1 $ and $ X = \PP^n $, namely $ h^1 (C,f^* T_{\tilde X}(-\epsilon)) $
for $ \epsilon \in \{ 0,1 \} $, where $ f^* T_{\tilde X}(-\epsilon) $ is to be
understood as $ (f^* T_{\tilde X}) \otimes \OO_C (-\epsilon) $.

\begin {lemma} \label {coh-1}
  Let $ C = \PP^1 $, $ X = \PP^n $, $ f: C \to \tilde X $ a morphism, and
  $ \epsilon \in \{ 0,1 \} $.
  \begin {enumerate}
  \item If $ f(C) \not\subset E $ or $f$ is a constant map then $ h^1
    (C,f^* T_{\tilde X}(-\epsilon)) = 0 $.
  \item If $ f(C) \subset E $ and the map $ f: C \to E \cong \PP^{n-1} $ has
    degree $ e>0 $ then
      \[ h^1 (C,f^* T_{\tilde X}(-\epsilon)) = e-1+\epsilon. \]
  \end {enumerate}
\end {lemma}

\begin {proof}
  (i): If $f$ is a constant map, the assertion is trivial, so let us assume
  that the homology class of the map is $ f_*[C] = dH'-eE' $ with $ d>0 $, $ e
  \ge 0 $. We have $ d = f^* H $ and $ e = f^* E $, and since $ f^* (H-E) $ is
  an effective divisor on $C$, it follows that $ e \le d $.

  By lemma \ref {triangle}, there is a commutative diagram
  \[ \begin {array}{c}
    f^* \pi^* T_X (-e-\epsilon) \ra f^* \pi^* T_X (-\epsilon) \\
    ~~~~~~ \da \\
    ~~~~~~ f^* T_{\tilde X} (-\epsilon)
  \end {array} \]
  where, in particular, the map $ f^* \pi^* T_X (-e-\epsilon) \to f^*
  T_{\tilde X} (-\epsilon) $ is injective and an isomorphism on a dense open
  subset of $C$. Hence we have an exact sequence
    \[ 0 \to f^* \pi^* T_X (-e-\epsilon)
         \to f^* T_{\tilde X} (-\epsilon)
         \to Q
         \to 0 \]
  with some sheaf $Q$ on $C$ which has zero-dimensional support. Therefore, to
  prove the lemma, it suffices to show that $ h^1 (C, f^* \pi^* T_X
  (-e-\epsilon)) = 0 $. But this follows easily from the Euler sequence,
  pulled back to $C$ and twisted by $ \OO_C (-e-\epsilon) $:
    \[ 0 \to \OO_C (-e-\epsilon)
         \to (n+1) \, \OO_C (d-e-\epsilon)
         \to f^* \pi^* T_X (-e-\epsilon)
         \to 0, \]
  since $ d-e-\epsilon \ge -\epsilon \ge -1 $.

  In particular, we also see that $ h^1 (C, f^* \pi^* T_X (-\epsilon)) = 0 $,
  which will be needed in the proof of part (ii).

  (ii): We consider the normal sequence
    \[ 0 \to T_E \to i^* T_{\tilde X} \to N_{E/\tilde X} \to 0. \]
  As $ N_{E/\tilde X} = \OO_E (-1) $, pulling back to $C$ and twisting by
  $ \OO_C (-\epsilon) $ yields
    \[ 0 \to f^* T_E (-\epsilon)
         \to f^* T_{\tilde X} (-\epsilon)
         \to \OO_C (-e-\epsilon)
         \to 0. \]
  By the remark at the end of part (i), applied to $ E \cong \PP^{n-1} $
  instead of $ X = \PP^n $, we see that $ h^1 (C, f^* T_E (-\epsilon)) = 0 $.
  Since $ h^1 (C, \OO_C (-e-\epsilon)) = e-1+\epsilon $, the result follows.
\end {proof}

Finally, we consider the case where $C$ is a genus 0 prestable curve in the
sense of \cite {BM}, i.e.\ a curve with at most ordinary double points as
singularities and whose arithmetic genus is zero.

\begin {prop} \label {coh-2}
  Let $C$ be a genus 0 prestable curve, $ X = \PP^n $, and $ f: C \to \tilde X
  $ a morphism. Let $ e' = e'(C) $ be "the sum of the exceptional degrees of
  all components of $C$ which are mapped into $E$", i.e.\
  \[ e' := \sum_{C'} \;\, \{ \mbox {
         $ e \; | \; C' $ is an irreducible component of $C$
         such that $ f_* [C'] = e \cdot E' $
       }\}. \]
  Then $ h^1 (C, f^* T_{\tilde X}) \le e' $, with strict inequality holding
  if $ e'>0 $.
\end {prop}

\begin {proof}
  The proof is by induction on the number of irreducible components of $C$. If
  $C$ itself is irreducible, the statement follows immediately from lemma
  \ref {coh-1}.

  Now let $C$ be reducible, so assume $ C = C_0 \cup C' $ where $ C' \cong
  \PP^1 $, $ C_0 \cap C' = \{Q\} $, and where $ C_0 $ is a prestable curve
  for which the induction hypothesis holds. Consider the exact sequences
  \begin {gather*}
    0 \to f^* T_{\tilde X}
      \to f_0^* T_{\tilde X} \oplus {f'}^* T_{\tilde X}
      \stackrel {\varphi}{\to} f_Q^* T_{\tilde X}
      \to 0 \\
    0 \to {f'}^* T_{\tilde X} (-Q)
      \to {f'}^* T_{\tilde X}
      \stackrel {\psi}{\to} f_Q^* T_{\tilde X}
      \to 0
  \end {gather*}
  where $ f_0 $, $ f' $, and $ f_Q $ denote the restrictions of $f$ to $ C_0 $,
  $ C' $, and $Q$, respectively.

  From these sequences we deduce that
  \begin {gather*}
    \dim \coker H^0 (\varphi)
      = h^1 (C, f^* T_{\tilde X}) - h^1 (C_0, f_0^* T_{\tilde X})
         - h^1 (C', {f'}^* T_{\tilde X}) \\
    \dim \coker H^0 (\psi)
      = h^1 (C', {f'}^* T_{\tilde X} (-Q)) - h^1 (C', {f'}^* T_{\tilde X}).
  \end {gather*}
  Since we certainly have $ \dim \coker H^0 (\varphi) \le \dim \coker H^0
  (\psi) $, we can combine these equations into the single inequality
    \[ h^1 (C, f^* T_{\tilde X})
         \le h^1 (C_0, f_0^* T_{\tilde X})
           + h^1 (C', {f'}^* T_{\tilde X} (-Q)). \]
  Now, by the induction hypothesis on $ f_0 $, we have $ h^1 (C_0, f_0^*
  T_{\tilde X}) \le e'(C_0) $ with strict inequality holding if $ e'(C_0) > 0
  $. On the other hand, we get $ h^1 (C', {f'}^* T_{\tilde X} (-Q)) \le e'(C')
  $ by lemma \ref {coh-1}. Since $ e'(C) = e'(C_0) + e'(C') $, the proposition
  follows by induction.
\end {proof}

%%%%%%%%%%%%%%%%%%%%%%%%%%%%%%%%%%%%%%%%%%%%%%%%%%%%%%%%%%%%%%%%%%%%%%%%%%%%%%%

\section {The morphism $ \phi : \mnb \to \mnx $} \label {image-map}

The proposition \ref {coh-2} shows us that, at least in the case $ X=\PP^n $,
problems with nonvanishing $ h^1 (C, f^* T_{\tilde X}) $ only arise if some
irreducible components of $C$ are mapped into the exceptional divisor. Since
these components are contracted by the map $ \pi: \tilde X \to X $, we are led 
to study the relation between stable maps in the blow-up $ \tilde X $ and their
image in $X$.

Let $\mtx$ denote the substack of $\mnx$ of those stable maps $ (C,x_1,\dots,
x_N,f) $ where $ (C,x_1,\dots,x_N) $ has a fixed topology which is encoded in
the graph $\tau$ as introduced in \cite {BM}. We will not need the details of
this encoding in this paper, all that will be important for us is that there is
a stratification of $\mnx$ by the $\mtx$ for all possible $\tau$, such that in
each stratum we have a fixed structure of the singularities of the curve $C$,
and the marked points lie in fixed components of $C$. Since we only consider
stable maps, there is only a finite number of possible graphs $\tau$ for given
$N$ and $\beta$.

Note that this is \emph {not} the stack $ \bar M_\tau (X,\beta) $ as defined in
\cite {BM}.

By the functorial properties of moduli spaces of stable maps \cite {BM}, the
map $ \pi: \tilde X \to X $ induces morphisms $ \phi: \mnb \to \mnx $ for all
$e$, where $ \beta \in H_2^+ (X) $ and where we use the decomposition $ H_2
(\tilde X) = H_2 (X) \oplus \ZZ \cdot E' $. We may restrict these maps to the
case where the underlying curves have topology $\tau$, so we also get morphisms
  \[ \phi_\tau: \mtb \to \mnx. \]
As we have seen in the example above, the dimension of the stack $ \mnb $ may
be larger than its virtual dimension, which is
\begin {align*}
  \vdim \mnb &= \chi (C, f^* T_{\tilde X}) + N - 3 \\
             &= n - K_X \cdot \beta - e \, (n-1) + N - 3 \\
             &= \vdim \mnx - e \, (n-1).
\end {align*}
The aim of this section is to prove the following proposition:

\begin {prop} \label {dim-im}
  Let $ X=\PP^n $ and $ \phi: \mnb \to \mnx $ be the morphism as above. Then
    \[ \dim \phi (\mnb) \le \vdim \mnb. \]
  Moreover, if $ R \subset \mnb $ denotes the subspace of all reducible stable
  maps in $\mnb$ (i.e.\ the maps $ (C,x_1,\dots,x_N,f) $ with $C$ reducible),
  then
    \[ \dim \phi (R) < \vdim \mnb. \]
\end {prop}

\begin {proof}
  Since the moduli spaces $ \mtb $ form a stratification of $ \mnb $, it is
  enough to show the statement for the restricted maps $ \phi_\tau $. Hence we
  fix a topology $ \tau $ and associate to it the following numerical
  invariants:
  \begin {itemize}
   \item Let $S$ be the number of nodes of a curve with topology $\tau$. We
     divide this number into $ S = S_{EE} + S_{XX} + S_{XE} $, where $ S_{EE} $
     (resp.\ $ S_{XX} $, $ S_{XE} $) denotes the number of nodes joining two
     exceptional components of $C$ (resp.\ two non-exceptional components, or
     one exceptional with one non-exceptional component). Here and in the
     following we call an irreducible component of $C$ exceptional if it is
     mapped by $f$ into the exceptional divisor and it is not contracted by
     $f$, and non-exceptional otherwise.
   \item Let $P$ be the (minimal) number of additional marked points which are
     necessary to stabilize $C$. We divide the number $P$ into $ P = P_E + P_X
     $, where $ P_E $ (resp.\ $ P_X $) is the number of marked points that have
     to be added on exceptional components (resp.\ non-exceptional components)
     of $C$ to stabilize $C$.
   \end {itemize}
  Now let $ \cc = (C,x_1,\dots,x_N,f) \in \mtb $ be a stable map of topology
  $ \tau $, and let
    \[ T_{\cc} \phi_\tau: T_{\mtb,\cc}
                    \to T_{\mnx,\phi (\cc)} \]
  be the differential of the map $ \phi_\tau $ at the point $ \cc $.
  As we always work over the ground field of complex numbers, to prove the
  proposition it suffices to show that
    \[ \dim \im T_{\cc} \phi_\tau \le \vdim \mnb \]
  for all $ \cc $, and that strict inequality holds if $ \tau $ is a topology
  corresponding to reducible stable maps.

  The tangent space $ T_{\mtb,\cc} $ is given by the hypercohomology group
  \cite {K}
    \[ T_{\mtb,\cc} = \HH^1 (T'_C \to f^* T_{\tilde X}) \]
  where $ T'_C = T_C (- x_1 - \dots - x_N) $ and where we put the sheaves
  $ T'_C $ and $ f^* T_{\tilde X} $ in degrees 0 and 1, respectively. This
  means that there is an exact sequence
    \[ 0 \to H^0 (C,T'_C)
         \to H^0 (C,f^* T_{\tilde X})
         \to T_{\mtb,\cc}
         \to H^1 (C,T'_C) \]
  (note that the first map is injective because $ \cc $ is a stable map). By
  Riemann-Roch we get $ \chi (C,T'_C) = S + 3 - N $. Moreover, by proposition
  \ref {coh-2} we have
    \[ \dim H^0 (C, f^* T_{\tilde X}) \le \chi (C, f^* T_{\tilde X}) + e'
       \qquad (*) \]
  where $e'$ is the "sum of the exceptional degrees of the components of $C$"
  as introduced there. It follows that
  \begin {align*}
    \dim T_{\mtb,\cc} &\le \chi (C, f^* T_{\tilde X}) + e' + N - S - 3 \\
                      &= \vdim \mnb + e' - S.
  \end {align*}
  We will now study the map
    \[ \phi': H^0 (C,f^* T_{\tilde X}) / H^0 (C,T'_C)
              \to T_{\mnx,\phi(\cc)} \]
  induced by the composition of the maps $ H^0 (C,f^* T_{\tilde X}) \to
  T_{\mtb,\cc} \to T_{\mnx,\phi(\cc)} $ considered above. To prove the
  proposition, we will show that
    \[ \dim \ker \phi' \ge e' - S \]
  and that strict inequality holds in certain cases. Obviously, we may
  assume that $ e'-S \ge 0 $.

  Let $ C_0 $ be a maximal connected subscheme of $C$ consisting only of
  exceptional components of $C$. Let $ f_0 $ be the restriction of $f$ to
  $ C_0 $ and let $ Q_1,\dots,Q_a $ be the nodes of $C$ which join $C_0$ with
  the rest of $C$ (they are of type $ S_{XE} $). Now every section of $ f_0^*
  T_E (-Q_1- \cdots - Q_a) $ can be extended by zero to a section of $ f^*
  T_{\tilde X} $ which is mapped to zero by $ \phi' $ since these deformations
  take place entirely within the exceptional divisor. As $ E \cong \PP^{n-1} $
  is a convex variety, we have
    \[ h^0 (C_0, f_0^* T_E) = \chi (C_0, f_0^* T_E)
       = n-1 + n \cdot e'(C_0) \]
  and therefore we can estimate the dimension of the space of deformations
  that we have just found:
    \[ h^0 (C_0, f_0^* T_E (-Q_1- \cdots -Q_a))
       \ge n-1 + n \cdot e'(C_0) - (n-1) \, a. \]
  (The right hand side of this inequality may well be negative, but
  nevertheless the statement is correct also in this case, of course.)

  We will now add up these numbers for all possible $ C_0 $, say there are
  $B$ of them. The sum of the $ e'(C_0) $ will then give $ e' = e'(C) $, and
  the sum of the $a$ will give $ S_{XE} $. Note that there is a $ P_E
  $-dimensional space of infinitesimal automorphisms of $C$, i.e.\ a subspace
  of $ H^0 (C,T'_C) $, included in the deformations that we have just found,
  and that these do not give non-trivial elements in the kernel of $ \phi' $.
  Therefore we have
  \begin {align*}
    \dim \ker \phi'
      &\ge B \, (n-1) + n e' - (n-1) S_{XE} - P_E \\
      &=   (n-1) \cdot (B + e' - S_{XE}) + e' - P_E \\
      &\ge B + e' - S_{XE} + e' - P_E
         \qquad \qquad \mbox {($ B+e'-S_{XE}\ge 0 $ since $ e' \ge S $)} \\
      &=   e' - S + (B + e' + S_{EE} - P_E) + S_{XX}. \qquad (+)
  \end {align*}
  Hence, it is certainly sufficient to show that $ P_E \le B + e' + S_{EE} $.

  To show this, we look at the exceptional components of $C$ where marked
  points have to be added to stabilize $C$. We have to distinguish three cases:
  \begin {itemize}
    \item Components on which two points have to be added, and whose (only)
      node is of type $ S_{EE} $: those give a contribution of 2 to $ P_E $,
      but they also give at least 1 to $e'$ and to $ S_{EE} $ (and every
      node of type $ S_{EE} $ "belongs" to at most one such component).
    \item Components on which two points have to be added, and whose (only)
      node is of type $ S_{XE} $: those give a contribution of 2 to $ P_E $,
      but they also give at least 1 to $e'$ and to $B$ (since such a component
      alone is one of the $ C_0 $ considered above).
    \item Components on which only one point has to be added: those give a
      contribution of 1 to $ P_E $, but they also give at least 1 to $e'$.
  \end {itemize}
  This finishes the proof of the "$\le$"-part of the lemma.

  To show the strict inequality if $C$ is reducible, we again distinguish two
  cases:
  \begin {itemize}
    \item $ e' > 0 $: Then, by lemma \ref {coh-2}, the inequality $(*)$ is
      already strict.
    \item $ e' = 0 $: If $C$ is reducible and $ e'=0 $, then we must have
      $ S_{XX} > 0 $, hence we get the strict inequality by $(+)$.
  \end {itemize}
  This completes the proof.
\end {proof}

%%%%%%%%%%%%%%%%%%%%%%%%%%%%%%%%%%%%%%%%%%%%%%%%%%%%%%%%%%%%%%%%%%%%%%%%%%%%%%%

\section {Geometric interpretation of the Gromov-Witten invariants on
  $ \tilde X $} \label {geometric}

We start this section by proving some moving lemmas which will be needed to
show that, in favourable cases, the intersection product on the moduli space
which defines the Gromov-Witten invariants can be made transverse.

\begin {lemma} \label {bertini-0}
  Let $X$ be a scheme of finite type and $ f: X \to \PP^m $ a morphism. Then,
  for a generic hyperplane $ H \subset \PP^m $, we have:
  \begin {enumerate}
    \item $ f^{-1}(H) $ is (empty or) of pure codimension 1 in $X$.
    \item If $X$ is smooth then the divisor $ f^{-1}(H) $ is a smooth
      subscheme of $X$ counted with multiplicity one.
  \end {enumerate}
\end {lemma}

\begin {proof}
  See \cite {J}, cor.\ 6.11.
\end {proof}

\begin {lemma} \label {bertini-1}
  Let $X$ be a scheme of finite type, $Y$ a smooth, connected, projective
  scheme, and $ f: X \to Y $ a morphism. Let $L$ be a base point free linear
  system on $Y$. Then, for generic $ D \in L $, we have:
  \begin {enumerate}
    \item $ f^{-1}(D) $ is (empty or) purely 1--codimensional.
    \item If $X$ is smooth then the divisor $ f^{-1}(D) $ is a smooth
      subscheme of $X$ counted with multiplicity one.
  \end {enumerate}
\end {lemma}

\begin {proof}
  The base point free linear system $L$ on $Y$ gives rise to a morphism
  $ s: Y \to \PP^m $ where $ m = \dim L $. Composing with $f$ yields a
  morphism $ X \to \PP^m $, and the divisors $ D \in L $ correspond to the
  inverse images under $s$ of the hyperplanes in $ \PP^m $. Hence, the
  statement follows from lemma \ref {bertini-0}, applied to the map $ X \to
  \PP^m $.
\end {proof}

\begin {lemma} \label {bertini-2}
  Let $X$ be a Deligne-Mumford stack of finite type, $ Y_i $ smooth,
  connected, projective schemes, and $ f_i : X \to Y_i $ morphisms for
  $ i = 1,\dots,N $. Let $ \alpha_i \in A^{c_i}(Y_i) $ be cycles of
  codimensions $ c_i \ge 1 $ on $ Y_i $ that can be written as intersection
  products of divisors on $ Y_i $
    \[ \alpha_i = [D'_{i,1}] \cdot \;\, \cdots \;\, \cdot [D'_{i,c_i}]
         \qquad \mbox {($ i = 1,\dots,N $)} \]
  such that the complete linear systems $ | D'_{i,j} | $ are base point free
  (this always applies, for example, in the case $ Y_i = \PP^n $). Let $ c =
  c_1 + \cdots + c_N $. Then, for almost all $ D_{i,j} \in | D'_{i,j} | $, 
  we have:
  \begin {enumerate}
    \item $ V_i := D_{i,1} \cap \cdots \cap D_{i,c_i} $ is smooth of
      pure codimension $ c_i $ in $ Y_i $, and the intersection is transverse.
      In particular, $ [V_i] = \alpha_i $.
    \item $ V := f_1^{-1}(V_1) \cap \cdots \cap f_N^{-1}(V_N) $ is of pure
      codimension $c$ in $X$. In particular, if $ \dim X < c $ then $ V =
      \emptyset $.
    \item If $ \dim X = c $ and $X$ contains a dense, open, smooth substack
      $U$ such that each geometric point of $U$ has no nontrivial
      automorphisms then $V$ consists of exactly $ (f_1^*\alpha_1 \cdots
      f_N^*\alpha_N)[X] $ points of $X$ (which lie in $U$ and are counted with
      multiplicity one).
  \end {enumerate}
\end {lemma}

\begin {proof}
  (i) follows immediately by recursive application of lemma \ref {bertini-0}
  to the schemes $ Y_i $.

  If $X$ is a scheme, then (ii) follows by recursive application of lemma \ref
  {bertini-1}. If $X$ is a Deligne-Mumford stack, then it has an \'etale cover
  $ S \to X $ by a scheme $S$, so (ii) holds for the composed maps $ S \to X
  \to Y_i $. But since the map $ S \to X $ is \'etale, the statement is also
  true for the maps $ X \to Y_i $.

  A Deligne-Mumford stack $U$ whose generic geometric point has no nontrivial
  automorphisms always has a dense open substack $U'$ which is a scheme (see
  e.g.\ \cite {V}. To be more precise, $U$ is a functor and hence an
  algebraic space (\cite {DM}, ex.\ 4.9), but an algebraic space always
  contains a dense open subset $U'$ which is a scheme (\cite {Kn}, p.\ 25)).
  Since $U'$ is dense in $X$ and therefore has smaller dimension, applying
  (ii) to the restrictions $ f_i |_{X \backslash U'} : X \backslash U' \to Y_i
  $ gives that $V$ is contained in the smooth scheme $U'$, hence it suffices
  to consider the restrictions $ f_i |_{U'} : U' \to Y_i $. But since $U'$
  is a smooth scheme, we can apply lemma \ref {bertini-1} (ii) recursively and
  get the desired result.
\end {proof}

We are now ready to give a geometric interpretation of the Gromov-Witten
invariants on $ \tilde \PP^n $ with only non-exceptional classes. First, we
show that the invariants can be thought of as certain numbers of stable maps
to $ \tilde \PP^n $.

\begin {prop} \label {finite}
  Let $ X=\PP^n $, $ 0 \neq \beta \in H_2^+ (X) $ and $ e \in \ZZ $. Let $
  \alpha_1,\dots,\alpha_N \in A^{\ge 1}(X) $ such that $ \sum_i \codim
  \alpha_i = \vdim \mnb $. Choose generic subschemes $ V_i \subset X $ with
  $ [V_i] = \alpha_i $ in the sense of lemma \ref {bertini-2}, which do not
  meet the blown up point $P$, such that the $ V_i $ may also be regarded as
  subschemes of $ \tilde X $.

  Then the number of stable maps $ (C,x_1,\dots,x_N,f) \in \mnb $ passing
  through the $ V_i $ (i.e.\ such that $ f(x_i) \in V_i $) is finite and equal
  to the Gromov-Witten invariant
    \[ \gwinv {\tilde X}{\pi^* \alpha_1,\dots,\pi^* \alpha_N}{\beta - eE'}. \]
  Moreover, all these stable maps are contained in the space $ RM'(\tilde X,
  \beta - eE') $ (see definition \ref {relevant}), and they are all counted
  with multiplicity one in the Gromov-Witten invariant.
\end {prop}

\begin {remark}
  In particular, if $ e<0 $ and $ \beta \neq 0 $ it follows that
    \[ \gwinv {\tilde X}{\pi^* \alpha_1,\dots,\pi^* \alpha_N}{\beta - eE'}
         = 0 \]
  for all $ \alpha_1,\dots,\alpha_N \in A^*(X) $, because there are no
  irreducible curves in $ \tilde X $ with homology class $ \beta - eE' $ for
  $ e<0 $, $ \beta \neq 0 $.
\end {remark}

\begin {proof}
  For $ i=1,\dots,N $ consider the commutative diagram
  \[ \begin {CD}
         \mnb        @>{\phi}>>    \mnx   \\
    @V{\tilde p_i}VV            @VV{p_i}V \\
       \tilde X      @>{\pi}>>       X
  \end {CD} \]
  where $ p_i $, $ \tilde p_i $ are the evaluation maps at the $i$-th marked
  point. By definition \cite {B}, the Gromov-Witten invariant mentioned in the
  proposition is equal to
  \[ \gwinv {\tilde X}{\pi^* \alpha_1,\dots,\pi^* \alpha_N}{\beta - eE'}
       = (\tilde p_1^* \pi^* \alpha_1 \cdots \tilde p_N^* \pi^* \alpha_N)
           \cdot [\mnb]^{virt} \qquad (*) \]
  where $ [\mnb]^{virt} $ denotes the virtual fundamental class. As shown in
  proposition \ref {coh-2}, we have $ h^1 (C, f^* T_{\tilde X}) = 0 $ whenever
  $C$ is irreducible. Therefore, on this part of the moduli space, the virtual
  fundamental class coincides with the usual one (since, by construction of
  virtual fundamental classes, this can be checked locally). This means that
  if $ I \subset \mnb $ is the substack consisting of irreducible maps and
  $ \bar I $ its closure, then one can write
    \[ [\mnb]^{virt} = [\bar I] + \gamma \qquad \in A_{\vdim \mnb} (\mnb) \]
  where $ \gamma $ is some cycle in $ \mnb $ whose support is entirely
  contained in the substack $ R \subset \mnb $ of reducible stable maps.
  But, by proposition \ref {dim-im}, we have $ \dim \phi(R) < \vdim \mnb $,
  therefore, by lemma \ref {bertini-2} (ii) applied to the restrictions $ p_i
  |_{\phi(R)}: \phi(R) \to X $ one concludes that, for generic choice of the
  $ V_i $, there are no stable maps in $ \phi(R) $ passing through the $ V_i
  $. This means that there are no reducible stable curves in $ \mnb $ passing
  through the $ V_i $. In particular, the contribution in $ (*) $ coming from
  the cycle $ \gamma $ vanishes, and we have
    \[ \gwinv {\tilde X}{\pi^* \alpha_1,\dots,\pi^* \alpha_N}{\beta - eE'}
         = (\tilde p_1^* \pi^* \alpha_1 \cdots \tilde p_N^* \pi^* \alpha_N)
           \cdot [\bar I]. \]
  So we only have to evaluate an intersection product with the usual
  fundamental class on $ \bar I $. By lemma \ref {bertini-2} (ii) this means
  that the number of stable maps passing through the $ V_i $ is finite and we
  are simply counting the number of such maps, although we do not yet know
  whether they are counted with multiplicity one.

  We do know, however, that all stable maps passing through the $ V_i $ are
  irreducible, which means, for example, that we can restrict ourselves to the
  case $ e \ge 0 $ since otherwise there certainly are no such curves. But we
  can restrict this even further. For example, we can assume that none of
  these maps is a finite covering map: for each irreducible finite covering
  map $ f: \PP^1 \to \tilde X $ of degree $ a>1 $ and homology class $ \beta
  - eE' $, there is also an irreducible stable map $ f': \PP^1 \to \tilde X $
  of homology class $ (\beta - eE')/a $. But the moduli space of maps of
  homology class $ (\beta - eE')/a $ is smaller than that of $ \beta - eE' $,
  since, if $ \beta = dH' $ for some $ d>0 $,
  \begin {align*}
    \vdim \mnb &- \vdim \mnbb {(\beta - eE')/a} \\
      &= - K_X \cdot \beta - e (n-1)
         + \frac 1 a \; (K_X \cdot \beta + e \, (n-1)) \\
      &= (d \, (n+1) - e \, (n-1)) \cdot (1 - \frac 1 a) \\
      &= [ (d-e) \cdot (n+1) + 2e ] \cdot (1 - \frac 1 a) \\
      &> 0 \qquad \mbox {since $ d \ge e \ge 0 $}
  \end {align*}
  (note that on the space of irreducible curves, the virtual dimension
  coincides with the actual one). Hence by lemma \ref {bertini-2} (ii), for
  generic $ V_i $ there are no maps $ f' $ as above passing through the
  $ V_i $, and therefore there are also no finite covering maps $ f: \PP^1
  \to X $.

  So, in the terminology of definition \ref {relevant}, we are only counting
  irreducible relevant maps. Let us denote the subspace of irreducible relevant
  maps by $ Z \subset \mnb $. Certainly, the locus $ RM'_N (\tilde X, \beta -
  eE') \subset Z $ of the maps $ (C,x_1,\dots,x_N,f) $ in $Z$ where $ f^{-1}
  (f(x_i)) = \{ x_i \} $ for all $i$ is dense in $Z$. Since $ RM'_N (\tilde X,
  \beta - eE') $ is smooth and, by lemma \ref {rel-2}, the geometric points in
  $ RM'_N (\tilde X, \beta - eE') $ have no non-trivial automorphisms, we can
  apply lemma \ref {bertini-2} (iii) to the restricted maps $ p_i|_Z: Z \to
  \tilde X $ and the statement of the proposition follows.
\end {proof}

Finally, we use the results of section \ref {map-image} to reformulate this
proposition in terms of embedded curves in $ \PP^n $ with multiple points,
which is our main result.

\begin {prop} \label {meaning}
  Let $ d>0 $, $ e \ge 0 $ and $ \alpha_1,\dots,\alpha_N \in A^{\ge 1}
  (\PP^n) $ such that
    \[ \sum_i (\codim \alpha_i - 1 ) = d \, (n+1) - e \, (n-1) + n - 3. \]
  Let $ P \in \PP^n $ be a point. Choose generic subschemes $ V_i \subset
  \PP^n $ with $ [V_i] = \alpha_i $ in the sense of lemma \ref {bertini-2}.

  Then the number of rational curves in $ \PP^n $ (purely 1-dimensional
  subvarieties birational to $ \PP^1 $) of degree $d$ which have non-empty
  intersection with all $ V_i $ and which pass through the point $P$ with
  total multiplicity $e$, where each such curve $C$ is counted with
  multiplicity
    \[ \sharp (C \cap V_1) \cdots \sharp (C \cap V_N), \]
 is equal to the Gromov-Witten invariant on $ \tilde \PP^n $
    \[ \gwinv {\tilde \PP^n}{\pi^*\alpha_1,\dots,\pi^*\alpha_N}{dH'-eE'}. \]
\end {prop}

\begin {proof}
  By proposition \ref {finite}, for generic $ V_i $ the Gromov-Witten
  invariant mentioned in the proposition counts maps $ (C,x_1,\dots,x_N,f) $
  in $ RM'_N (\tilde \PP^n, dH'-eE') $ with multiplicity one which pass
  through the $ V_i $, and there are no further stable maps in $ \mnbb
  {dH'-eE'} $ passing through the $ V_i $. By proposition \ref {rel-1} and
  lemma \ref {rel-2}, this can be reformulated by saying that we are counting
  irreducible relevant curves in $ \tilde \PP^n $ with multiplicity one which
  pass through the $ V_i $, i.e.\ rational curves $ C \subset \tilde \PP^n $
  of homology class $ dH'-eE' $ together with points $ p_i \in C $ such that $
  p_i \in V_i $ for all $i$. This is obviously the same as counting rational
  curves $ C \subset \tilde \PP^n $ meeting all the $ V_i $, and counting them
  with multiplicity $ \sharp (C \cap V_1) \cdots \sharp (C \cap V_N) $.

  But, using the strict transform of curves $ C \subset \PP^n $ under the
  blow-up $ \pi: \tilde \PP^n \to \PP^n $, there is a one-to-one correspondence
  between rational curves in $ \PP^n $ of degree $d$ which pass through $P$
  with multiplicity $e$ and rational curves in $ \tilde \PP^n $ with homology
  class $ dH'-eE' $. Of course, the property of meeting the subvarieties $ V_i
  $ is not affected by this strict transform, hence the proposition follows.
\end {proof}

%%%%%%%%%%%%%%%%%%%%%%%%%%%%%%%%%%%%%%%%%%%%%%%%%%%%%%%%%%%%%%%%%%%%%%%%%%%%%%%

\section {Calculation of the Gromov-Witten invariants of $ \tilde \PP^n $}
  \label {calculation}

To compute the Gromov-Witten invariants of $ \tilde \PP^n $, we recall the
First Reconstruction Theorem of Kontsevich and Manin.

\begin {prop} \label {reconstruction}
  Let $X$ be a smooth projective variety such that the algebraic part of $
  H^*(X) $ is generated as a ring by divisor classes on $X$. Then all
  Gromov-Witten invariants $ \gwinv {X}{\alpha_1,\dots,\alpha_N}{\beta} $
  with $ \alpha_i \in A^* (X) $ and $ \beta \in H_2^+ (X) $ can be
  reconstructed from those with $ N=3 $ and $ \alpha_3 \in A^1 (X) $.
\end {prop}

\begin {proof}
  \cite {KM} theorem 3.1, applied to the Gromov-Witten classes constructed in
  \cite {B}. An explicit algoritm to compute the invariants is also given in
  \cite {KM}.
\end {proof}

Now we apply this result to the case of $ X = \tilde \PP^n $. We set $ H_i =
H^i $, $ E_i = -(-E)^i $, and choose
  \[ B = \{ H_0, H=H_1, H_2,\dots,H_n,E=E_1,E_2,\dots,E_{n-1} \} \]
as a basis of $ A^* (\tilde \PP^n) $. In the following, we consider only
invariants whose classes are in this basis.

\begin {prop} \label {calc}
  All Gromov-Witten invariants on $ \tilde \PP^n $ can be computed
  recursively by proposition \ref {reconstruction} using the initial data
  \begin {enumerate}
   \item $ \gwinv {\tilde \PP^n}{pt,pt,H}{H'} = 1 $, where $ pt $ denotes
     the class of a point,
   \item $ \gwinv {\tilde \PP^n}{E_{n-1},E_{n-1},E}{E'} = -1 $,
   \item $ \gwinv {\tilde \PP^n}{\alpha_1,\alpha_2,\alpha_3}{H'-E'} = 1 $
     if $ \alpha_i \in B $ and $ \codim \alpha_1 + \codim \alpha_2 + \codim
     \alpha_3 = n+2 $,
   \item $ \gwinv {\tilde \PP^n}{\alpha_1,\alpha_2,\alpha_3}{dH'+eE'} = 0 $ in
     all other cases where $ \codim \alpha_3 = 1 $ and $ \alpha_i \in B $.
  \end {enumerate}
\end {prop}

\begin {proof}
  Let $ \gwinv {\tilde \PP^n}{\alpha_1,\alpha_2,\alpha_3}{dH'+eE'} $ be an
  invariant with $ \codim \alpha_3 = 1 $. Let $ i = \codim \alpha_1 $ and
  $ j = \codim \alpha_2 $, such that $ 1 \le i,j \le n $.

  Case 1: $ d=0 $, $ e>0 $. Then the dimension condition reads
    \[ i+j = e \, (n-1) + n - 1 = (e-1)(n-1) + 2n - 2. \]
  Since these curves are contained in the exceptional divisor, the
  Gromov-Witten invariant is zero if there is a point class (or any other
  non-exceptional class) among the $ \alpha_i $. Hence we may assume that $ i
  + j \le 2n-2 $. But then it follows that $ e=1 $ and $ \alpha_1 = \alpha_2
  = E_{n-1} $, and we are in case (ii).

  To prove (ii), note that for maps $ f: C \to E $ of degree 1 into the
  exceptional divisor, we have $ h^1 (C, f^* T_{\tilde \PP^n}) = 0 $ by
  proposition \ref {coh-2}, hence the corresponding moduli stack is smooth of
  the expected dimension, and its virtual fundamental class coincides with the
  usual one. Now consider the invariant
    \[ \gwinv {\tilde \PP^n}{H_{n-1}-E_{n-1},H_{n-1}-E_{n-1},H-E}{E'}. \]
  The classes $ H_{n-1} - E_{n-1} $ and $ H-E $ are represented on $ \tilde
  \PP^n $ by the strict transform of a line (resp.\ hyperplane) in $ \PP^n $
  passing through $P$. These intersect the exceptional divisor transversally in
  a point (resp.\ in a hyperplane in $E$), hence this Gromov-Witten invariant
  simply counts the number of lines in $E$ through two points in $E$ (and
  intersecting a hyperplane in $E$), which is 1.

  Note that Gromov-Witten invariants $ \gwinv {\tilde \PP^n}{\alpha_1,\dots,
  \alpha_N}{eE'} $ with $ e>0 $ vanish if one of the $ \alpha_i $ is a
  non-exceptional class, since one can choose a subvariety $ V_i \subset \PP^n
  $ representing $ \alpha_i $ which does not pass through the exceptional
  divisor, such that there are no stable maps of homology class $ eE' $ passing
  through $ V_i $. Hence, by linearity of the Gromov-Witten invariants it
  follows that
    \[ \gwinv {\tilde \PP^n}{E_{n-1},E_{n-1},E}{E'}
         = - \gwinv {\tilde \PP^n}{H_{n-1}-E_{n-1},H_{n-1}-E_{n-1},H-E}{E'}
         = - 1, \]
  which proves (ii).

  Case 2: $ d>0 $. Then we must have
    \[ 0 \le (H-E)(dH'+eE') = d+e \]
  and the dimension condition is
  \begin {align*}
    i+j &= d \, (n+1) + e \, (n-1) + n - 1 \\
        &= (d+e)(n-1) + 2d + n - 1. \qquad (*)
  \end {align*}
  If $ d+e > 0 $, then we have $ i+j \ge 2d + 2n - 2 $, therefore it follows
  that $ d=1 $ and $ \alpha_1 = \alpha_2 = pt $, which is (i). The statement
  of (i) follows immediately from proposition \ref {meaning}.

  Now assume that $ e=-d $. Since any curve in $ \PP^n $ of degree $d$ passing
  with multiplicity $d$ through a given point is a union of lines, the image of
  any stable map of homology class $ dH' - dE' $ is simply a line passing
  through the exceptional divisor (note that it is not possible to have a union
  of some lines since such a curve would not be connected in $ \tilde \PP^n $).

  So we consider the moduli space of lines in $ \tilde \PP^n $ intersecting
  $E$, which is canonically isomorphic to $ E \cong \PP^{n-1} $ itself. The
  condition that such a line meets a generic linear subspace of codimension $k$
  in $ \PP^n $ or $E$ is given by a codimension $ k-1 $ linear subspace in the
  moduli space $E$. Therefore, the moduli space of lines intersecting $E$
  and two linear subspaces corresponding to $ \alpha_1 $ and $ \alpha_2 $ is
    \[ \dim E - (i-1) - (j-1) = n - i - j + 1. \qquad (+) \]
  So to get a non-zero invariant, we must have $ i + j \le n+1 $. Therefore
  from $(*)$ it follows that $ d=1 $, in which case $(+)$ is zero, such that
  there is exactly one line satisfying the incidence conditions. Note that by
  lemma \ref {coh-2} the moduli stack of stable maps of homology class $ H'-E'
  $ is smooth of the expected dimension, and we indeed count the number of
  lines intersecting $E$ and two classes representing $ \alpha_1 $ and $
  \alpha_2 $. This proves (iii).

  Finally, we have also shown that (i)--(iii) are the only non-zero invariants
  on $ \PP^n $ with $ N=3 $ and $ \codim \alpha_3 = 2 $, which proves (iv).
\end {proof}

A few numerical examples of invariants can be found in tables \ref {invP2}
to \ref {invexc}. There are some remarks that can be made about these
numbers:
\begin {itemize}
 \item Let $ \alpha_1,\dots,\alpha_k $ be non-exceptional classes, i.e.\
   $ \alpha_i \in \pi^* A^* (\PP^n) $. If $ d>0 $ then
     \[ \gwinv {\tilde \PP^n}{\alpha_1,\dots,\alpha_n}{dH'-E'}
          = \gwinv {\tilde \PP^n}{\alpha_1,\dots,\alpha_n,pt}{dH'}. \]
   This follows, for example, from proposition \ref {meaning}, since both
   sides count the number of rational curves in $ \PP^n $ of degree $d$
   which meet subvarieties representing the $ \alpha_i $ and one additional
   point.

   This is no longer true if one of the $ \alpha_i $ is exceptional.
 \item Similarly, if $ \alpha_1,\dots,\alpha_k $ are non-exceptional classes,
   $ d>0 $ and $ e<0 $, then
     \[ \gwinv {\tilde \PP^n}{\alpha_1,\dots,\alpha_n}{dH'-eE'} = 0. \]
   Again, this follows from proposition \ref {meaning} and is no longer
   true if one of the $ \alpha_i $ is exceptional. Indeed, the invariants need
   not even be positive, so that there can be no geometric interpretation of
   these invariants as numbers of curves (at least not in an obvious way).
 \item Again let $ \alpha_1,\dots,\alpha_k $ be non-exceptional classes. If
   $ e=0 $ then
     \[ \gwinv {\tilde \PP^n}{\alpha_1,\dots,\alpha_n}{dH'}
          = \gwinv {\PP^n}{\alpha_1,\dots,\alpha_n}{dH'}. \]
   This can also be deduced from proposition \ref {meaning}. The Gromov-Witten
   invariants of $ \PP^2 $ can be found e.g.\ in \cite {KM} 5.2.1., and those
   of $ \PP^3 $ and $ \PP^4 $ in \cite {JS}.
 \item The numbers of rational curves of degree $d$ in $ \PP^2 $ passing
   through $ 3d-3 $ points in the plane and in addition through the point $P$
   with multiplicity two, which we have computed as Gromov-Witten invariants
   on $ \tilde \PP^2 $ (case $e=2$ in table \ref {invP2}), have already been
   calculated by R. Pandharipande in \cite {P} by different methods, and indeed
   the numbers agree.
 \item Consider the invariants $ \gwinv {\tilde \PP^2}{pt,\dots,pt}{dH'-
   (d-1)E'} $ for $ d>1 $, where we have $ 2d $ point classes in the invariant.
   A curve $C$ of degree $d$ in $ \PP^2 $ passing with multiplicity $ d-1 $
   through a point $P$ has genus
     \[ \frac 1 2 (d-1) (d-2) - \frac 1 2 (d-1)(d-2) = 0, \]
   i.e.\ it is always a rational curve. Hence the space of degree $d$ rational
   curves with a $(d-1)$-fold point in $P$ is simply a linear system of the
   expected dimension. Therefore we have
     \[ \gwinv {\tilde \PP^2}{pt,\dots,pt}{dH'- (d-1)E'} = 1 \]
   which can also be seen in table \ref {invP2}.
 \item From the dimension of the moduli space of stable maps of homology class
   $ dH' - eE' $ in $ \PP^2 $, it can be seen that the dimension of the space
   of degree $d$ rational curves in $ \PP^2 $ which pass through a given point
   with multiplicity $e$ is $ 3d-1-e $. This can also be understood
   geometrically by a dimension count as follows: if $ C \subset \PP^2 $ is a
   curve of degree $d$ having ordinary $k_i$-fold points in the points $ P_i $
   and no further singularities, the genus $g$ of the normalization of $C$ is
     \[ g = \frac 1 2 (d-1)(d-2) - \sum_i \frac 1 2 k_i \, (k_i-1) \]
   Hence, if we already have an $e$-fold point at $P$, to get a rational curve
   we need additional
     \[ \frac 1 2 (d-1)(d-2) - \frac 1 2 e \, (e-1) \]
   double points. Now consider the linear system of degree $d$ curves in
   $ \PP^2 $, which has dimension $ \frac 1 2 (d+1)(d+2) - 1 $. Having an
   $e$-fold point at $P$ gives $ \frac 1 2 e \, (e+1) $ conditions, and
   each of the above additional double points imposes one more constraint
   (since the points where these double points occur are not specified). Hence
   the dimension of the space of curves we considered is
     \[ \frac 1 2 (d+1)(d+2) - 1 - \frac 1 2 e \, (e+1)
          - \left[ \frac 1 2 (d-1)(d-2)
                     - \frac 1 2 e \, (e-1) \right]
        = 3d-1-e. \]
\end {itemize}

%%%%%%%%%%%%%%%%%%%%%%%%%%%%%%%%%%%%%%%%%%%%%%%%%%%%%%%%%%%%%%%%%%%%%%%%%%%%%%%

% Some Gromov-Witten invariants of $ \tilde \PP^2 $ and $ \tilde \PP^3 $

\begin {table}[p]
  \[ \begin {array}{|l|r|r|r|r|r|r|r|} \hline
          & d=1 & d=2 & d=3 & d=4 &  d=5  &   d=6    &     d=7     \\ \hline
    e = 0 &   1 &   1 &  12 & 620 & 87304 & 26312976 & 14616808192 \\ \hline
    e = 1 &   1 &   1 &  12 & 620 & 87304 & 26312976 & 14616808192 \\ \hline
    e = 2 &   0 &   0 &   1 &  96 & 18132 &  6506400 &  4059366000 \\ \hline
    e = 3 &   0 &   0 &   0 &   1 &   640 &   401172 &   347987200 \\ \hline
    e = 4 &   0 &   0 &   0 &   0 &     1 &     3840 &     7492040 \\ \hline
    e = 5 &   0 &   0 &   0 &   0 &     0 &        1 &       21504 \\ \hline
    e = 6 &   0 &   0 &   0 &   0 &     0 &        0 &           1 \\ \hline
  \end {array} \]
  \caption {Some Gromov-Witten invariants $ \gwinv {\tilde \PP^2}{pt,\dots,
    pt}{dH'-eE'} $, where the number of point classes in the invariant is $ 3d-
    1 - e $. By proposition \ref {meaning}, these are the numbers of degree $d$
    rational curves in $ \PP^2 $ meeting $ 3d-1-e $ generic points in the
    plane, and in addition passing through the point $P$ with multiplicity
    $e$.}
  \label {invP2}
\end {table}

\begin {table}[p]
  \[ \begin {array}{|l|r|r|r|r|r|r|r|r|} \hline
          & d=1 & d=2 & d=3 & d=4 & d=5 & d=6  &  d=7   &   d=8   \\ \hline
    e = 0 &   1 &   0 &   1 &   4 & 105 & 2576 & 122129 & 7397760 \\ \hline
    e = 1 &   1 &   0 &   1 &   4 & 105 & 2576 & 122129 & 7397760 \\ \hline
    e = 2 &   0 &   0 &   0 &   0 &  12 &  384 &  23892 & 1666128 \\ \hline
    e = 3 &   0 &   0 &   0 &   0 &   0 &    0 &    620 &   72528 \\ \hline
    e = 4 &   0 &   0 &   0 &   0 &   0 &    0 &      0 &       0 \\ \hline
  \end {array} \]
  \caption {Some Gromov-Witten invariants $ \gwinv {\tilde \PP^3}{pt,\dots,
    pt}{dH'-eE'} $, where the number of point classes in the invariant is $ 2d-
    e $. These are the numbers of rational curves of degree $d$ in $ \PP^3 $
    meeting $ 2d-e $ generic points, and in addition passing through the point
    $P$ with multiplicity $e$.}
  \label {invP3}
\end {table}

\begin {table}[p]
  \[ \begin {array}{|l|r|r|r|r|} \hline
          &  d=1   &    d=2    &      d=3      &       d=4       \\ \hline
    e =-3 &   2925 &   4849635 &   25767926176 & 362956315020486 \\ \hline
    e =-2 &    -68 &    -35832 &     -89070592 &   -730861150688 \\ \hline
    e =-1 &      3 &       342 &        382720 &      1793900214 \\ \hline
    e = 0 &      0 &         0 &         -2332 &        -5810112 \\ \hline
    e = 1 &      1 &        -3 &            40 &           23825 \\ \hline
    e = 2 &      0 &         0 &             4 &             960 \\ \hline
    e = 3 &      0 &         0 &             0 &              45 \\ \hline
  \end {array} \]
  \caption {Some Gromov-Witten invariants $ \gwinv {\tilde \PP^3}{E_2,\dots,
    E_2}{dH'-eE'} $, where $ E_2 = -E^2 $ and the number of classes in the
    invariant is $ 4d-2e $.}
  \label {invexc}
\end {table}

\end {document}